\documentclass{article}

\newcommand\aw{{{\mathbf A}_W}}
\newcommand\np{{{\cal N}{\cal P}}}
\newcommand\nn{\{0,1\}^*}

\begin{document}

\title{Expansions of pseudofinite structures and circuit and proof complexity}
\author{Jan Kraj\'{\i}\v{c}ek}

\date{Faculty of Mathematics and Physics\\
Charles University in Prague}

\maketitle

{\em 
I am honored to have a chance
to contribute to this volume celebrating a personal anniversary
of Albert Visser.}

\begin{abstract}
I shall describe a general model-theoretic task to construct
expansions of pseudofinite structures and discuss several examples of particular 
relevance to computational complexity. Then I will
present one specific situation where finding a suitable expansion would imply that, 
assuming a one-way permutation exists, the computational class $\np$ is not closed
under complementation.
\end{abstract}

Consider the following situation: ${\mathbf M}$ is a nonstandard model of true arithmetic
(in the usual language of arithmetic
$0, 1, +, \cdot, \le$), $n$ is a nonstandard element of ${\mathbf M}$,
$L$ is a finite language and $W \in {\mathbf M}$ is its 
interpretation on the universe $[n]=\{1,\dots ,n\}$; $W$ can be identified with 
a subset of 
$[n^k]$ for some $k \in {\mathbf N}$. We shall denote the resulting structure
$\aw$; it is coded by an element
of ${\mathbf M}$ that is $\le 2^{n^k}$.
Without a loss of generality we shall assume that $L$ contains 
constants $1,n$, the ordering relation $\le$ interpreted as in ${\mathbf M}$, and ternary relation symbols
$\oplus$ and $\odot$ for the graphs of addition and multiplication inherited
from ${\mathbf M}$. Because ${\mathbf M}$ is a model of true arithmetic
$\aw$ is pseudofinite: it satisfies the theory of all finite $L$-structures.

Paris and Dimitracopoulos \cite{ParDim82} studied the problem of
for how large $m > n$
does the theory of the arithmetic structure on $[n]$ determine the theory
of the arithmetic structure on $[m]$ and proved that it does not for $m=2^n$.
They also pointed out various links between questions of this type 
and complexity theory problems around the collapse of the polynomial time
hierarchy. Ajtai \cite{Ajt83} showed (among other similar results)
that if $M$ is a countable nonstandard model of
PA and $L$ is finite then for any
$L$-structure $\aw$ there are two sets $U, U' \subseteq [n]$, both elements of 
${\mathbf M}$,
such that ${\mathbf M}$ thinks that $|U|$ is odd and $|U'|$ is even while,
as structures, $(\aw,U) \cong (\aw,U')$ (the isomorphism is not in ${\mathbf M}$, of
course)\footnote{We shall discuss another example with parity later.}.
Kraj\'{\i}\v cek and Pudl\' ak \cite{KP89} showed that for any nonstandard $t\le n \in 
{\mathbf M}$ one can construct ${\mathbf M}' \supseteq {\mathbf M}$ containing
a proof of contradiction
in PA of length bounded by $n^t$ without adding any new
elements to interval $[0,n]$. M\' at\' e \cite{Mate} considered 
the full second order
structure on $[n]$ (coded in ${\mathbf M}$) and reformulated the 
statement that $\np \neq co\np$ as a statement about non-preservation
of the theory of the structure
in an expansion coming from a model ${\mathbf M}' \supseteq {\mathbf M}$.

The most interesting results of this kind (to this author) were
obtained initially 
by Ajtai \cite{Ajt83,Ajt88}. In the first paper he established
that parity of $n$ bits cannot be computed by $AC^0$ circuits (proved
independently by Furst, Saxe and Sipser \cite{FSS}) and he reports 
there that his first proof of the lower bound was by model theory
of arithmetic although he eventually chose 
to present the result combinatorially. In the second paper he proved that propositional formulas $PHP_m$ 
formalizing the pigeonhole
principle do not have polynomial size constant depth Frege proofs.
That proof is by constructing a suitable model of arithmetic\footnote{The original
manuscript was only about models of bounded arithmetic. After
Ajtai learned that Paris and Wilkie \cite{ParWil85} linked provability of PHP
in bounded arithmetic with a conjecture of Cook and Reckhow \cite{CR}
that formulas $PHP_m$ are hard for Frege systems,
shown eventually false by Buss \cite{Bus87}, he added a few hand-written pages showing how his result implies a lower bound for constant depth Frege systems.}.
We shall
discuss these two examples in the next section.

In this, mostly expository, note
we are interested in the general question of how to construct 
expansions of $\aw$ with particular
properties. Before formulating this more specifically
we will consider in the next section three examples.
The examples go back to Ajtai \cite{Ajt83,Ajt88} and (essentially)
M\' at\' e \cite{Mate}
but two of them are not presented in the literature with enough details
and are formulated with unnecessarily strong hypotheses.
In the subsequent section we shall discuss a specific open problem 
whose solution would have
interesting implications for computational complexity.

The note is self-contained modulo a basic knowledge of logic and complexity theory.
Notions not explained here can be likely found in \cite{kniha}.

\section{Examples}

\noindent
{\bf Parity example.}
It is well-known that the parity of a string of bits
cannot be computed by $AC^0$ circuits, \cite{Ajt83,FSS}.
That is, 

\begin{enumerate}

\item [(1)] For any $d \geq 1$ and large enough $m$, any depth $d$, size $\le m^d$
circuit with $m$ inputs must compute erroneously the parity of some $m$-bit string.

\end{enumerate}
Computability by $AC^0$ circuits is equivalent to first-order definability
in the presence of an extra structure (of a fixed signature, the extra
structure depending just on the size
of the universe). In particular, (1) is equivalent to
\begin{enumerate}

\item [(2)] For any finite language $L$ and any formula $\Phi(X)$
in the language $L(X)$, $L$ augmented by a unary predicate $X(x)$, for $m$ 
large enough and any $L$-structure ${\mathbf B}$ with universe $[m]$
it holds:

\begin{itemize}

\item There is $U \subseteq [m]$ such that the equivalence:
$$
({\mathbf B}, U) \models \Phi(U)\ \mbox{ iff }\ \mbox{$|U|$ is odd}
$$
fails.

\end{itemize}
\end{enumerate}  
Using overspill in ${\mathbf M}$, (2) is equivalent to
the same statement for any $L$-structure on any nonstandard $[n]$ in ${\mathbf M}$
(with $U \in {\mathbf M}$ and its parity
defined in ${\mathbf M}$). And that can be further
formulated as follows.
For $u \in [n]$ denote $U^{<u}:=\{v \in U\ |\ v < u\}$ and
$U^{\le u}:=\{v \in U\ |\ v \le u\}$. 

\begin{enumerate}

\item [(3)] For any nonstandard $n \in {\mathbf M}$, any finite language $L$,
any $L(X)$ formula $\Phi(X)$ and any $L$-structure $\aw \in {\mathbf M}$ with universe $[n]$
it holds:

There is $U \subseteq [n]$, $U \in {\mathbf M}$, and $u \in U$
such that the following holds:
\begin{itemize}

\item [(a)]
$$
(\aw, U^{<u}) \models \Phi(U^{<u})\ \mbox{ iff }\ 
(\aw, U^{\le u}) \models \Phi(U^{\le u})\ .
$$
\item [(b)]
For $t \in U^{<u}$,
$$
(\aw, U^{<t}) \models \Phi(U^{<t})
\ \mbox{ iff }\ 
(\aw, U^{\le t}) \models \neg \Phi(U^{\le t})\ .
$$
\end{itemize}
\end{enumerate}
Let $Y(x)$ be a new unary predicate
and consider first-order
theory $T_1^\Phi$ in the language $L(X,Y)$, $L(X)$ augmented by $Y$, with axioms:
\begin{enumerate}

\item The least number principle axioms:
$$
\exists x\ \alpha(x) \rightarrow 
(\exists x \forall y < x\ \alpha(x) \wedge \neg \alpha(y))
$$
for all formulas $\alpha$ in the language $L(X,Y)$
(as we will evaluate formulas over a structure with
universe $[n]$ they are de facto bounded),

\item $\exists x \in X,\ (x \in Y) \not \equiv \Phi(X^{\le x})$,

\item axiom $\Psi(Y)$, where $\Psi$ is the formula:
$$
Y \subseteq X \wedge \min X \in Y \wedge
(\forall y\in X,\ suc_X(\min X)=y \rightarrow y \notin Y)\ \wedge\ \ \ \ \ \ \ \ \ \ 
$$
$$
\ \ \ \ \ \ \ \ \ \ \ \ \ \ \ \ \forall x, y \in X,\ suc_X(suc_X(x))=y \rightarrow
x \in Y \equiv y \in Y
$$
where for $x \in X$, $suc_X(x)$ is $\min \{z \in X\ |\ x < z\}$, if it exists.

\end{enumerate}

\medskip
\noindent
{\bf Claim 1:} {\em Statement (3) for a
given language $L$, formula $\Phi$, set $U \in {\mathbf M}$ and some $u \in U$
is equivalent to
the existence of $V \subseteq [n]$ ($V$ not necessarily in ${\mathbf M}$)
such that the expanded structure $(\aw,U,V)$ satisfies $T_1^\Phi$
($U$ interprets $X$ and $V$ interprets $Y$).}

\medskip

For the only if direction note that for $U$ and $u \in U$ satisfying (3) 
a suitable $V$ can be defined in ${\mathbf M}$
already: take for $V$ the subset of $U$ consisting of its elements on odd
positions. Axiom 1. holds because ${\mathbf M}$ satisfies it 
for all formulas, axiom 3. holds by the definition of $V$ and axiom
2. holds as it is witnessed by $x:=u$.

For the if direction assume that 
$(\aw,U,V)\models T_1^\Phi$. Take for $u$ the minimal $x$ witnessing
axiom 2.; it exists by the least number principle.
Utilizing axiom 3. we see 
that the pair $U, u$ satisfies statement (3).

\bigskip
\noindent
{\bf PHP example.}
In this example we aim at a proof complexity lower bound. Given $m \geq 2$,
consider a propositional formula $PHP_m$ formed using atoms $p_{i j}$, $i \in [m]$
and $j \in [m-1]$ that is the disjunction of the following formulas:
\begin{itemize}

\item $\bigvee_i \bigwedge_j \neg p_{i j}$,

\item $\bigvee_{i_1 \neq i_2} \bigvee_j (p_{i_1 j} \wedge p_{i_2 j})$,

\item $\bigvee_{i} \bigvee_{j_1 \neq j_2} (p_{i j_1} \wedge p_{i j_2})$.

\end{itemize}
Having a falsifying assignment $p_{i j}:=a_{i j} \in \{0,1\}$ we could define
the graph of an injective map from $[m]$ to $[m-1]$:
$$
\{(i,j)\ |\ a_{i j} = 1\}
$$
which is impossible. So $PHP_m$ is a tautology.

The depth of a formula (in DeMorgan language) is the maximal number of
connected blocks of alike  
connectives on a path of subformulas of the formula (the depth of $PHP_m$
is thus $3$). A Frege system $F$ is a sound and implicationally complete finite
collection of inference rules and axiom schemes. Its depth $d$ subsystem $F_d$ is allowed
to use only formulas of depth $\le d$.

A path in a depth $d$ formula of size (i.e. the number of symbols) $s$ can be
naturally determined by a $d$-tuple of numbers $\le s$, giving $d$ pointers to which
subformula the path moves when alternating from one connective to another. 
This implies that such a formula can be
coded by a $d$-ary function on $[s]$ giving information about atoms (or constants) 
in which individual paths
end. Hence an $F_d$ proof of $PHP_n$ of size polynomial in 
$n$ can be coded by a relation $S$ on $[n]$. The qualification {\em naturally} 
used above means
that, given $d$, one can define in $([n],R,S)$ the satisfaction relation 
between  formulas in the proof $S$ and a truth assignment\footnote{Details left 
out in this example
can be found in \cite{kniha}.}$R$.

Using this, and overspill as before, the statement

\begin{enumerate}

\item [(4)]
For all $d \geq 3$ and all large enough $m \geq 2$
there is no $F_d$ proof of $PHP_m$ of size $\le m^d$,

\end{enumerate}
is equivalent to the following statement about ${\mathbf M}$ and any
nonstandard $n \in {\mathbf M}$:

\begin{enumerate}

\item [(5)]
for all $d \geq 3$ and any relation $S$ on $[n]$, $S \in {\mathbf M}$,
$S$ does not encode an
$F_d$ proof of $PHP_n$.

\end{enumerate}
Let $L$ be a finite language
and let $Z$ be a variable for a binary 
relation. Consider the theory $T_2$:
\begin{itemize}

\item Induction axioms for all $L(Z)$-formulas as in $T_1^\Phi$,

\item $\neg PHP(Z)$ axiom:
$$
[\forall x \exists y < n\ Z(x,y)]\ \wedge\ 
[\forall x < x' \forall y < n\ \neg Z(x,y) \vee \neg Z(x', y)]\ \wedge
$$
$$
\ \ \ \ \ \ \ \ \ \ \ 
[\forall x \forall y < y' < n\ \neg Z(x,y) \vee \neg Z(x, y')]\ .
$$

\end{itemize}

\medskip
\noindent
{\bf Claim 2:} {\em Assume that for any finite $L$ 
and any $L$-structure $\aw$ with universe $[n]$ (any nonstandard $n \in {\mathbf M}$) 
there exists an expansion
$(\aw,R)$ satisfying $T_2$ with $Z:=R$. Then statement (5) is true.}

\medskip

To see this assume that (5) fails, i.e. some relation $S$ on $[n]$
encodes an $F_d$-proof of $PHP_n$. Let $(\aw,R)$ be
a model of $T_2$, with $W$ containing $S$. Using the truth definition
for depth $\le d$ formulas show that under the assignment
$$
p_{i j}:=1\ \mbox{ if $R(i,j)$ and }
p_{i j}:=0\ \mbox{ otherwise}
$$
all propositional axioms and the formula $\neg PHP_n$ are true.
Hence, by induction, there has to be an inference whose all
hypotheses are true while its concussion is not. But that is impossible.

Note that Claim 2 is formulated as an implication and not equivalence
as Claim 1; we shall return to this issue in the next section.

\bigskip
\noindent
{\bf TAUT example.}
Let $TAUT \subseteq \nn$ be the set of propositional tautologies in the
DeMorgan language.
A propositional proof system (abbreviated to PPS) in the sense of Cook and
Reckhow \cite{CR} is a polynomial time binary relation $P$ on $\nn$ such that
\begin{itemize}
\item $\forall x,y\ P(x,y) \rightarrow x \in TAUT$ (soundness),

\item $\forall x\ \in TAUT \rightarrow \exists y\ P(x,y)$ 
(completeness).
\end{itemize}
A PPS $P$ is p-bounded if there exists $k \geq 1$ such that the $\exists y$
in the completeness can be bounded by $|y|\le |x|^k$. A p-bounded PPS
exists iff $\np = co\np$, see \cite{CR}.
The main task is therefore to establish for all PPSs a super-polynomial lengths-of-proofs
lower bound. This may be far away at present
but lower bounds for specific PPSs have
interesting consequences as well (e.g. for independence results 
or for SAT algorithms).

As before we shall consider strings of length polynomial in $n$ (in particular, formulas)
coded by relations on $[n]$. For a given PPS
$P$ and $k \geq 1$ there exists by Fagin's theorem\footnote{Invoking just Fagin's theorem here is
not enough for various properties one needs often from $\Theta$. The formalization needs to be
"natural" again. The claim below holds for arbitrary $\Theta$ defining the relation, though.} a formula 
$\Theta(X,Y,Y')$ such that $\exists Y'\Theta(X,Y,Y')$ defines
on $[n]$ the relation $P(x,y) \wedge |y|\le |x|^k$.

Denote by $SAT(Z,X)$ a first order formula defining
the satisfaction relation between an assignment $Z$ and a DNF formula $X$
(we want to avoid coding of evaluations of general formulas in this discussion). 
Define a theory
$T_3$ with axioms:
\begin{itemize}
\item $\neg SAT(Z, X)$,

\item $\forall X',Y,Y',Z'\ \Theta(X',Y,Y') \rightarrow SAT(Z', X')$.
\end{itemize}
Note that it is a $\Pi^1_1$ theory, not first order as $T_1^\Phi$ and $T_2$.

\medskip
\noindent
{\bf Claim 3:}{\em 
Assume that for any finite $L$, any $L$-structure $\aw$ with universe $[n]$ and any standard $k \geq 1$
there is $F \subseteq [n]$, $F \in {\mathbf M}$ a DNF formula that is a tautology in ${\mathbf M}$,
such that $(\aw,F)$ has an expansion $(\aw,F,R)$
satisfying $T_3$ with $X:= F$ and $Z:=R$. Then $P$ is not
p-bounded.
}

\medskip

If $P$ were p-bounded with exponent $|x|^k$ then by overspill
the formula $F$ would have a $P$-proof in ${\mathbf M}$ of size polynomial in $n$,
and $\Theta(F,S,S')$ would hold for some relations $S, S'\in {\mathbf M}$ on $[n]$
that we can put into $W$. Hence an expansion that is a model of $T_3$
for $X:=F$ is impossible.

\section{Discussion}

Claims 2 and 3 can be established as equivalence statements using the
theory of propositional translations of $\Pi^1_1$ theories (cf.\cite[Chpt.9]{kniha}).
For that argument to work one does not need that ${\mathbf M}$ is
a model of true arithmetic but only that
\begin{itemize}

\item [(a)] $2^{n^t}$ exists in ${\mathbf M}$ for some nonstandard $t \le n$,

\item [(b)] ${\mathbf M}$ satisfies bounded arithmetic theory $R^1_2$ (which yields
induction for $\Sigma^1_1$ formulas on $\aw$),

\item [(c)] ${\mathbf M}$ is countable.

\end{itemize}
Ajtai \cite{Ajt07,Ajt11} formulated a general existence theorem for theories
$T$ going well beyond first order or $\Pi^1_1$ theories. Such $T$ can be not only second order
or third order, etc., but it can be a finite set theory over Ur-elements $[n]$
(and even more), and he allows not only expansions of $\aw$ but expansions
of end-extensions of $\aw$. The existence of such a model of $T$ is characterized
by the non-existence of a proof of contradiction in $T$ that is - in a specific,
rather technical, sense - definable over $\aw$. We refer the reader to
Garl\'{\i}k \cite{Gar} who found a simpler and more conceptual proof
of Ajtai's theorem. The construction needs ${\mathbf M}$ satisfying (a)-(c) 
above and also
\begin{itemize}

\item [(d)] $L$ is finite.

\end{itemize}
In Claims 1 - 3 we stipulated that $L$ is finite in order to avoid the
discussion how it is coded. In these claims $L$ can be, in fact, infinite
as long as $\aw$ is coded in ${\mathbf M}$.
But in \cite{Ajt07,Ajt11,Gar} the
hypothesis (d) is needed.

The intended goal of Claims 1 - 3 is to 
offer a strategy how to prove a lower bound in complexity
theory by constructing a suitable expansion of $\aw$. This has been done
by Ajtai \cite{Ajt83,Ajt88} for the lower bounds explained in 
Claims 1 and 2. Ajtai \cite{Ajt83,Ajt88}
works in a model of PA but the construction needs only
assumptions (a) and (c) above and a variant of (b):
\begin{itemize}

\item [(b')] ${\mathbf M}$ is a model of the theory PV and of the weak 
pigeonhole principle for p-time functions (as represented by PV-terms), denoted
WPHP(PV).

\end{itemize}
See \cite[Sec.15.2]{kniha} for how the WPHP(PV) is used.

\bigskip

It is a challenge to construct a suitable model of $T_3$ in 
the situation of Claim 3 for strong PPS $P$ (or even for all PPSs). 
In an attempt to meet the challenge particular models ${\mathbf M'}$
(extending a cut in ${\mathbf M}$)
were constructed in \cite{k2,Kra-hardcore} and \cite{Kra-nw} (two different constructions).
They were constructed under the assumptions that a one-way permutation exists.
These models ${\mathbf M'}$ satisfy the following conditions:

\begin{enumerate}

\item $2^{2^{n^\epsilon}}$ exists in ${\mathbf M'}$ for some standard $\epsilon > 0$
(this is stronger that (a)),

\item ${\mathbf M'}$ is a model of the true $\forall\Pi^b_1$ theory 
in the language of PV (this does not imply either (b) or (b')).
In particular, all PPSs are sound in ${\mathbf M'}$.

\item $L$ is infinite and coded in ${\mathbf M'}$ but $\aw$ is not coded in ${\mathbf M'}$,

\item there is a DNF formula $F \subseteq [n]$, $F \in {\mathbf M}$ (the original model of true
arithmetic) that has the form:
$$
\bigvee_{i \in [n]} \varphi_i(x, y^i)
$$
with each $\varphi_i(x, y^i)$ having the form
$$
\psi_i(x, y^i) \rightarrow \eta_i(y^i)
$$
and such that:
\begin{itemize}

\item [(i)] $F$ is a tautology in both ${\mathbf M}$ and ${\mathbf M'}$,

\item [(ii)] $x, y^1, \dots, y^n$ are mutually disjoint tuples of variables,

\item [(iii)] there exists assignments $A \in W$ and $B^i \in W$ for all $i \in [n]$
such that $\varphi_i(A,B^i)$ fails, i.e. $SAT( (A,B^i), \neg \varphi_i)$, 
and thus $SAT((A,B^i), \psi_i) \wedge SAT(B^i, \neg \eta_i)$,
hold in ${\mathbf M'}$.

\end{itemize}

\end{enumerate}
We now derive more properties of ${\mathbf M'}$ using the assumption that a PPS $P$
is p-bounded and sufficiently strong\footnote{We need that $P$ contains resolution,
$P$-proofs of true sentences can be constructed by a p-time algorithm
(hence its soundness is true in ${\mathbf M'}$), and that $P$ 
simulates modus ponens and substitutions of constants with only a polynomial increase
in the proof length.}. We shall use the specific form of the formula $F$ used in 
\cite{k2,Kra-nw}; we will formulate its properties explicitly but the interested reader is
assumed to learn the definition of the formula in \cite{k2,Kra-nw} and check that it has
the stated properties. 

\medskip
\noindent
{\bf Claim 4:}{\em
Assume $P$ is p-bounded and sufficiently strong. Then, in ${\mathbf M'}$, the following 
facts hold:
\begin{itemize}

\item [(iv)] There is a $P$-proof of $\bigvee_{i \in [n]} \varphi_i(A,y^i)$ in ${\mathbf M'}$
coded by a relation on $[n]$,

\item [(v)] for each $i \in [n]$ there are $P$-proofs of 
$(\psi_i(x,z) \wedge \psi_i(x, z')) \rightarrow z \equiv z'$ and of 
$\neg \varphi_i(A,y^i)$ 
in ${\mathbf M'}$ coded by relations on $[n]$,

\item [(vi)] formula $\bigwedge_{i \in [n]} \psi_i(A,y^i)$ is not $P$-refutable in ${\mathbf M'}$.

\end{itemize}
}
Statement (iv) follows from the p-boundedness of $P$ and (i) above. 
The first part of statement (v) follows again from the p-boundedness of $P$
and the fact that formulas 
$(\psi_i(x,z) \wedge \psi_i(x, z')) \rightarrow z \equiv z'$
are tautologies in ${\mathbf M}$. It follows that 
$\psi_i(A,B^i) \rightarrow 
\psi_i(A, y^i)$ and thus $\psi_i(A, y^i)$ 
(using that the true sentences $\psi_i(A, B^i)$ have $P$-proofs)
and $\neg \eta_i(B^i)$ have $P$-proofs in ${\mathbf M'}$ too.
Statement (vi) is valid because of the specific formulas $\psi_i$ have
the property that
it holds in ${\mathbf M}$ that
$$
\forall x \exists y=((y)_1, \dots, (y)_n) \bigwedge_{i \in [n]} \psi(x,(y)_i)\ .
$$
This implies that in ${\mathbf M}$ for no $a$ and $k \in {\mathbf N}$
there is a $P$-proof of 
$\neg \bigwedge_{i \in [n]} \psi_i(a,y^i)$ of size $\le n^k$, and by p-boundedness of $P$
these facts have $P$-proofs in ${\mathbf M}$ (and hence in ${\mathbf M'}$)
and $P$ is sound in ${\mathbf M'}$.

\bigskip

Ideally we would like to bring the existence of ${\mathbf M'}$ with properties 1.-4.(with (iv)-(vi) added)
to a contradiction. That would imply that the hypothesis of the existence of a one-way permutation
used in \cite{k2,Kra-nw} contradicts the hypothesis of p-boundedness of $P$ used in Claim 4, entailing a conditional
lower bound for (possibly all) $P$.  In light of Claim 4 it seems natural to try to extend model ${\mathbf M'}$ 
by adding a satisfying assignment for $\bigwedge_i \psi(A, y^i)$. 

The existence of a satisfying assignment for $\bigwedge_i \psi(A, y^i)$ would follow if
${\mathbf M'}$ would satisfy a collection scheme for $\Sigma^1_1$ formulas on $\aw$. But 
that is unlikely as the argument of Cook and Thapen \cite{CT} implies that the true $\forall\Pi^b_1$
theory does not prove the scheme (unless factoring is not hard), and the same argument 
implies that one cannot argue just using the $\forall\Pi^b_1$
theory of ${\mathbf M'}$ that it has an extension where the collection holds
(we would also need that it adds no new elements into $[n]$).
The construction of ${\mathbf M'}$ in \cite{Kra-nw} depends just on the $\forall\Pi^b_1$
theory of ${\mathbf M}$ and hence it is unlikely that it can be used again with the
ground model being ${\mathbf M'}$. 

The use of the property that we add no new elements into $[n]$ in Claims 1 - 3 was
solely to preserve the first order theory of $\aw$ (in the preceding paragraph it would also
imply that (iii) above remains true in an extension of ${\mathbf M'}$). 
There is an alternative for arranging that.
In the forcing from \cite{k2,Kra-hardcore} one naturally adds many new elements into $[n]$
(and, in fact, does not include all of $[n]$ from 
${\mathbf M}$ into ${\mathbf M'}$)
but the first order theory of $\aw$  is nevertheless often preserved.
The set-up of the method presupposes some approximate counting available in 
the ground model (here ${\mathbf M'}$) and the model is assumed to be
$\aleph_1$-saturated. The former can be arranged as one
modify the constructions so that the so called dual WPHP  
for p-time functions, dWPHP(PV), is true in ${\mathbf M'}$ too and that
yields some approximate counting by Je\v r\' abek \cite{Jer1,Jer2}.
As for the latter condition: some saturation property of ${\mathbf M'}$ could be
arranged (cf. \cite{Kra-satur}) if one could modify the forcing construction
of ${\mathbf M' }$ so that it is defined by a compact family (in the sense of
\cite{k2}) of random variables. That would be possible if one could
establish a hard-core lemma for the computation model underlying the construction. 
That, together with the fact (see 3. above) that $\aw$ is not coded
in ${\mathbf M'}$, seem to be the main technical obstacles to apply
Claim 3 to a general PPS.

\bigskip
\noindent
{\bf Ackowledgements:} I thank Neil Thapen (Prague) for comments on
a draft of this paper.

\bigskip
\noindent
{\bf Mailing address:}

Department of Algebra

Faculty of Mathematics and Physics

Charles University

Sokolovsk\' a 83, Prague 8, CZ - 186 75

The Czech Republic

{\tt krajicek@karlin.mff.cuni.cz}


\begin{thebibliography}{44}

\bibitem {Ajt83}
M.~Ajtai, $\Sigma^1_1$ - formulae on finite structures,
{\em Annals of Pure and Applied Logic}, {\bf 24},
(1983), pp.1-48.

\bibitem {Ajt88} 
M.~Ajtai, The complexity of the pigeonhole principle, 
in: {\em Proc. IEEE 29$^{\mbox{th}}$ Annual Symp. on
Foundation of Computer Science}, (1988), pp. 346-355.

\bibitem {Ajt07}
M.~Ajtai,
Generalizations of the Compactness Theorem and Godel's 
Completeness  Theorem  for  Nonstandard  Finite  Structures, 
in: Proc. of the
4th international conference on Theory and applications of 
models of computation, (2007), pp.13-33.

\bibitem {Ajt11}
M.~Ajtai,
A Generalization of Godel's Completeness Theorem for 
Nonstandard Finite Structures, manuscript (2011).


\bibitem {Bus87}
S.~R.~Buss,
The propositional pigeonhole principle has polynomial size 
Frege proofs, {\em J. Symbolic Logic}, {\bf  52}, (1987), pp.916-927.

\bibitem {CR}
S.~A.~Cook, and Reckhow, The relative
efficiency of propositional proof systems, 
{\em J. Symbolic Logic}, {\bf 44(1)}, (1979), pp.36-50.

\bibitem {CT}
S.~A.~Cook and N.~Thapen,
The strength of replacement in weak arithmetic, {\em ACM Transactions on Computational Logic}, 
{\bf Vol 7:4}, (2006), pp.749-764.

\bibitem {FSS}
M.~Furst, J.~B.~Saxe, and M.~Sipser, 
Parity, circuits and the polynomial-time hierarchy, 
{\em Math. Systems Theory}, {\bf 17}, (1984), pp.13-27.

\bibitem {Gar}
M.~Garl\'{\i}k, 
A New Proof of Ajtai's Completeness Theorem for Nonstandard Finite Structures,
{\em Archive for Mathematical Logic}, {\bf 54(3-4)}, (2015), pp. 413-424.

\bibitem {Jer1}
E.~Je\v r\' abek,
Dual weak pigeonhole principle, Boolean complexity, and derandomization, 
{\em Annals of Pure and Applied Logic}, {\bf 129}, (2004), pp. 1–37.

\bibitem {Jer2}
E.~Je\v r\' abek,
Approximate counting in bounded arithmetic, {\em J. of Symbolic Logic},
{\bf 72(3)}, (2007), pp. 959–993.



\bibitem{kniha}
J.~Kraj\'{\i}\v cek, {\em Bounded arithmetic, propositional
logic, and complexity theory},  Encyclopedia of Mathematics
and Its Applications, Vol. {\bf 60}, Cambridge University Press,
(1995).

\bibitem{Kra-nw}
J.~Kraj\'{\i}\v cek, 
On the proof complexity of the Nisan-Wigderson generator based on a 
hard $NP \cap coNP$ function,
J. of Mathematical Logic, Vol.11 (1), (2011), pp.11-27. 


\bibitem{Kra-hardcore}
J.~Kraj\'{\i}\v cek, 
Pseudo-finite hard instances for a student-teacher game with a Nisan-Wigderson generator,
Logical methods in Computer Science, Vol. 8 (3:09) 2012, pp.1-8.


\bibitem{Kra-satur}
J.~Kraj\'{\i}\v cek, 
A saturation property of structures obtained by forcing with a compact family of random variables,
Archive for Mathematical Logic, 52(1), (2013), pp.19-28. 


\bibitem{k2}
J.~Kraj\'{\i}\v cek, 
{\em Forcing with random variables and proof complexity}, 
London Mathematical Society Lecture Note Series, No.382, Cambridge University Press, (2011).


\bibitem {KP89}
J.~Kraj\'{\i}\v cek and P. Pudl\'{a}k,
On the Structure of Initial Segments of Models of Arithmetic, 
{\em Archive for Mathematical Logic}, {\bf 28(2)}, (1989), pp, 91-98. 

\bibitem {Mate}
A.~M\' at\' e, 
Nondeterministic Polynomial-Time Computations and
Models of Arithmetic, {\em J. of ACM}, {\bf37(1)}, (1990),
pp.175-193.

\bibitem {ParDim82}
J.~Paris and C.~Dimitracopoulos, Truth definitions
for $\Delta_0$ formulae, in : {\em Logic and Algorithmic},
{\em l'Enseignement Mathematique}, {\bf 30}, (1982), pp.318-329,
Gen\'{e}ve.

\bibitem {ParWil85} 
J.~Paris and A.~J.~Wilkie, 
Counting problems in
bounded arithmetic, in: {\em Methods in Mathematical Logic}, 
LNM 1130, (1985), pp.317-340.Springer.

\end{thebibliography}
\end{document}